\documentclass[11 pt]{article}  
\usepackage[utf8]{inputenc}
\usepackage{amsmath}
\usepackage{amsfonts}
\usepackage{amssymb}
\usepackage{graphicx}
\usepackage{mathrsfs}
\usepackage{upref,amsthm,amsxtra,exscale}
\usepackage{stmaryrd}
\usepackage{cite}
\usepackage[colorlinks=true,urlcolor=blue,
citecolor=red,linkcolor=blue,linktocpage,pdfpagelabels,
bookmarksnumbered,bookmarksopen]{hyperref}

\usepackage{fullpage}

\usepackage{subcaption}
\usepackage{caption}
\usepackage{cleveref}

\usepackage{enumitem}

\newtheorem{theorem}{Theorem}[section]

\newtheorem{lemma}[theorem]{Lemma}
\newtheorem{proposition}[theorem]{Proposition}

\numberwithin{equation}{section}

\theoremstyle{definition}
\newtheorem{remark}[theorem]{Remark}

\usepackage{color}
\usepackage[dvipsnames]{xcolor}

\usepackage[T1]{fontenc}
\usepackage[utf8]{inputenc}
\usepackage{authblk}

\title{\textbf{Normalized Solutions to Schr\"{o}dinger Equations with General Nonlinearities in Bounded Domains via a Global Bifurcation Approach}}
\author[a,b]{Wei Ji\thanks{Email: jiwei2020@amss.ac.cn} }
\affil[a]{\small Academy of Mathematics and Systems Science, Chinese Academy of Sciences, Beijing 100190, P.R. China}
\affil[b]{\small{University of Chinese Academy of Sciences, Beijing 100049, P.R. China}}
  \date{}

\begin{document}

\maketitle

\abstract{We obtain the existence, nonexistence and multiplicity of positive solutions with prescribed mass for nonlinear Schr\"{o}dinger 
equations in bounded domains via a global bifurcation approach. 
The nonlinearities in this paper can be mass supercritical, critical, 
subcritical or some mixes of these cases, and the equation can be autonomous or non-autonomous. 
This generalizes a result in Noris, Tavares and Verzini [\emph{Anal. PDE}, 7 (8) (2014) 1807-1838], where
the equation is autonomous with homogeneous nonlinearities. 
Besides, we have proven some orbital stability or instability results.}

\bigskip
\bigskip

\noindent\text{\textbf{Keywords:}} Global bifurcations, positive solutions, normalized solutions, orbital stability and instability.
\medskip

\noindent\text{\textbf{MSC2020:} 35A15, 35B09, 35B40, 35B35}


\medskip

\section{Introduction}

Suppose that $\Omega$ is a $C^1$ bounded domain in $\mathbb R^N$.
In this paper, we consider the positive solutions 
to the following semilinear elliptic problem: 
\begin{align}\label{c:intro}
  \begin{cases}-\Delta u +\lambda u= f(x,u)~~ \text{ in }\Omega,\\
    u=0, ~~ \text{ on }\partial \Omega,
    \end{cases}
\end{align}
with prescribed $L^2$-norm
\begin{equation}
    \int_{\Omega}\vert u \vert^2 =c. 
\end{equation}

In the last decades, significant progress has been made in the study of normalized solutions for the following Schr\"{o}dinger equation:
\begin{align}\label{c:introrn}
-\Delta u +\lambda u= f(x,u)\quad \quad u \in H^1(\mathbb R^N).
\end{align}
For example, see~\cite{normalizedinrn1}-\cite{normalizedinrn3}. Since the energy functional 
$$I(u)=\frac{1}{2}\int_{\mathbb R^N}|\nabla u|^2-\int_{\mathbb R^N}F(x,u)$$ restricted on 
$\{u \in H^1(\mathbb R^N):\|u\|^2_{L^2(\mathbb R^N)}=c\}$ is unbounded from
below when $f(x,u)$ is $L^2$-supercritical, most results are obtained from a natural constraint on the Pohozaev manifold 
(see \cite{normalizedinrn2}). 
Furthermore, the authors in~\cite{zhongxuexiu} considered the asymptotic behaviors
of the positive solutions of ~\eqref{c:introrn} as $\lambda \to 0^+$ and $\lambda \to +\infty$, and get a global branch
to study the positive normalized solutions. 

But the results of normalized solutions for Schr\"{o}dinger equations in bounded domains are still few. Indeed, the classical
natural constraint approach does not work since there is a boundary term that cannot be eliminated
in the Pohozaev indentity of the energy functional, restricted on $\{u \in H_0^1(\Omega):\|u\|^2_{L^2(\Omega)}=c\}$, 
of equation~\eqref{c:intro}. 
Moreover, since the bounded domains can not be scaled arbitrarily, there are some difficulties in studying the 
asymptotic behaviors of the positive solutions of~\eqref{c:intro}. 

Recently, by considering a dual-constraint problem, the authors in \cite{Noris1}-\cite{Noris4} studied the positive normalized solutions 
of autonomous Schr\"{o}dinger equations 
with homogeneous nonlinearity ($f(u)=|u|^{p-2}u$) in bounded domains. That is, the authors
focused on the variational problem
$$
\max\left\{\int_{\Omega}|u|^pdx:u\in H_0^1(\Omega), \int_{\Omega}u^2dx=1, \int_{\Omega}|\nabla u|^2dx=a\right\}
$$
and its asymptotic properties with respect to the parameter $a$. 
In fact, the method for proving the existence of solutions in these papers strictly depends on the uniqueness of the positive solution
and homogeneity of the nonlinearity of the equation. 
And L. J. Song in \cite{LSongCVPDE} considered the existence, nonexistence and multiplicity of the
positive normalized solution of Schr\"{o}dinger equations with the mass supcritical and Sobolev subcritical nonlinearities
relying on some bifurcation type arguments, under the condition that the equation admits only one positive solution.

In this paper, we present a new approach to constructing a global bifurcation of solutions to study the asymptotic behaviors
of the positive solutions of ~\eqref{c:intro} as $\lambda \to +\infty$, and then derive
 the existence, nonexistence and multiplicity of positive normalized solutions. 

\medskip

We assume on $\Omega$ and $f$:

\begin{itemize}

 \item[$(\Omega)$]$\Omega \in \mathbb R^N(N\geq 2)$ is convex in the $x_i$-direction for any $i=1,2,\ldots, N$ and $G$-invariant,  
 where $G$ is the group of isometries $G:=O(m_1)\times\ldots\times O(m_\ell)\subset O(N)$, with $\ell \geq 1, m_i\geq2$,
 for $i=1,\ldots,\ell$ such that $\sum_{i=1}^{\ell}m_i=N$, and $O(m_i)$ denotes the group of linear isometries of $\mathbb{R}^{m_i}$.

  \item[$(f_{1})$]$f(x,t) \in C^{1}(\Omega \times \mathbb R, \mathbb R)$,
  $f(x,0)=f_u(x,0)=0$ and there exist $\alpha,\beta$ satisfying $2<\alpha\leq\beta<2^*$
  where
  \begin{align}
  2^*=
\begin{cases}\frac{2N}{N-2}, \qquad  N \geq 3,
    \\ +\infty, \qquad  N=1,2,\end{cases}
  \end{align}
  
  such that
  \begin{align}\label{estimateA}
  0<(\alpha-1)f(x,u)u\leq f_u(x,u)u^2\leq(\beta-1)f(x,u)u\quad \text{a.e.} \quad x\in\Omega,u\neq0.
  \end{align}
  
  By \eqref{estimateA}, we know that
  \begin{align}\label{estimateB}
  0<\alpha F(x,u)\leq f(x,u)u\leq\beta F(x,u)\quad a.e.\quad x\in\Omega, u\neq0
  \end{align}
  where $F(x,u)=\int_0^uf(x,s)ds.$

  \item[$(f_{2})$] $f(x,1)$ and $f(x,-1)$ is bounded, 
  i.e. there exists a constant $M>0$ such that $|f(x,1)|<M$ and $|f(x,-1)|<M$, $\forall x \in \Omega$.

  \item[$(f_{3})$] $f(x,u)$ satisfies 
  \begin{align}\label{Q:intro3}
    f(x,u)= d(x)|u|^{p-2}u+g(x,u),
   \end{align}
   where $p \in (2,2^*)$,  $d(x)\geq0$, $d(0)>0$ , $ d(x) \in C^1(\bar \Omega) \cap L^\infty(\Omega)$, 
   and $g(x,u)$ satisfies
   \begin{align}\label{Q:intro3B}
    \lim_{t \to \infty}\frac{g(x,t)}{t^{p-1}}=0~~\text{and}~~\lim_{t \to 0}\frac{g(x,t)}{t}=0. 
   \end{align}
   Moreover, $f(x,\cdot)$ and $d(x)$ are nonincreasing in $|x_i|$ for all $i=1,\ldots,N$, and
   $f(hx,\cdot)=f(x,\cdot), d(hx)=d(x)$ for any $h\in G$, where $G$ is defined in $(\Omega)$.

 \end{itemize}

 By \eqref{Q:intro3} we know that $F(x,u)=\frac{1}{p} \int_{\Omega} d(x)|u|^{p}+G(x,u)$, where $G(x,u)=\int_0^u g(x,s)ds$.

 \begin{remark}\label{existoff}
  If $(f_1)$ and $(f_2)$ hold, it is standard to prove the existence of positive solutions of \eqref{c:intro} without normalized conditions via Mountain pass theorem (see \cite[Section 1.3]{Minimax}).
   \end{remark}

\medskip

Without loss of generality, we assume that $d(0)=1$.

For the equation
\begin{align}\label{global}
  -\Delta v + v&=|v|^{p-2}v \quad \text{ in }\mathbb R^N,
\end{align} where $p \in (2,2^*)$, 
\cite{Kwong} yields that \eqref{global} admits a unique positive (radial) solution $Q$. 
It is not difficult to know that $Q \in H_G^1(\mathbb R^N)$ and $Q$ is a ground state of \eqref{global},
where 
$$
H_G^1(\mathbb R^N)=\{u\in H^1(\mathbb R^N):gu=u, \forall g \in G\}, 
$$
and $G$ is defined in $(\Omega)$.

Set
$$
E:=\{(\lambda,u)\in\mathbb{R}\times H_0^1(\Omega)\setminus\{0\}:\int_\Omega(\nabla u\nabla\varphi-\lambda u\varphi)dx
=\int_\Omega f(x,u)\varphi dx,\forall\varphi\in H_0^1(\Omega)\}.
$$

Our main conclusions are as follows. 
 \begin{theorem}\label{Normalizedsolution}
  Let $N\geq 2$, $\lambda> -\lambda_1$, where $\lambda_1$ denotes the first Dirichlet eigenvalue of $-\Delta$ on $\Omega$,  
   $\Omega$ satisfy $(\Omega)$, and $f$ satisfy $(f_1)-(f_3)$. Assume that \eqref{c:intro} admits a
unique positive solution, then the following statements hold:
\begin{itemize}
  \item[$(1)$] 
  If $2+\frac{4}{N}<p<2^*$, then there exists some $b>0$ such that 
  \begin{itemize}
  \item[$(i)$] there exists $(\lambda, u_\lambda) \in E$ such that $u_\lambda>0$ and $\|u_\lambda\|_{L^2(\Omega)}^2=b$;
  \item[$(ii)$] for any $0<c<b$, there exist $(\lambda,u_\lambda),(\tilde{\lambda},u_{\tilde{\lambda}}) \in E$, 
  such that $u_\lambda>0$, $u_{\tilde{\lambda}}>0$ and 
  $\|u_\lambda\|_{L^2(\Omega)}^2=\|u_{\tilde{\lambda}}\|_{L^2(\Omega)}^2=c$;
  \item[$(iii)$]for any $c>b$, there exists no $(\lambda,u_{\lambda})\in E$ such that  
       $u_{\lambda}>0$ and $\|u_\lambda\|_{L^2(\Omega)}^2=c$.
    \end{itemize}

    \item[$(2)$]
    If $p=2+\frac{4}{N}$, then there exists $D \geq d=\|Q\|_{L^2(\mathbb R^N)}^2>0$ such that 
 \begin{itemize}
    \item[$(i)$] for any $0<c < d$, there exists $(\lambda,u_\lambda) \in E$, such that $u_\lambda>0$ and $\|u_\lambda\|_{L^2(\Omega)}^2=c$;
    \item[$(ii)$] for any $c>D$, there exists no $(\lambda,u_{\lambda})\in E$, such that $u_\lambda>0$ and $\|u_\lambda\|_{L^2(\Omega)}^2=c$;
  \end{itemize} 
  where $Q$ is the unique positive solution of \eqref{global}.
  
    \item[$(3)$]
    If $2<p<2+\frac{4}{N}$, then for any $0<c<+\infty$, there exists $(\lambda,u_\lambda) \in E$, 
    such that $u_\lambda>0$ and $\|u_\lambda\|_{L^2(\Omega)}^2=c$.

  \end{itemize}

\end{theorem}

\medskip
\begin{remark}  

  Regarding the hypotheses $2+\frac{4}{N}<p<2^*$ of $f$ in statement $(1)$ in Theorem \ref{Normalizedsolution}, 
  there are several possible cases: 
  \begin{itemize}
  \item[$(a)$] $f$ is $L^2$-supcritical, i.e. 
  \begin{align}\label{Q:intro4A}
  \lim_{t\to 0}\frac{f(|t|)}{|t|^{1+\frac{4}{N}}}=0~~~\text{and}~~~\lim_{t\to \infty}\frac{f(|t|)}{|t|^{2^*-1}}=0. 
  \end{align}
  \item[$(b)$] $f$ is mixed with $L^2$-supcritical and critical terms, such as $f=|u|^{p-2}u+|u|^{l_1-2}u$, 
  where $l_1=2+\frac{4}{N}<p<2^*$; or mixed with $L^2$-supcritical and subcritical terms, such as $f=|u|^{p-2}u+|u|^{l_2-2}u$, 
  where $2<l_2<2+\frac{4}{N}<p<2^*$.
 \end{itemize}
 Similar considerations apply to the hypotheses in statement $(2)$ as well. 
 \end{remark}  
 
   \medskip

To prove this result, we first study the following bifurcation phenomena.

\begin{theorem}\label{uatinfity}
  Let $N\geq 2$, $\lambda> -\lambda_1$, $\Omega$ satisfy $(\Omega)$, 
  and $f$ satisfy $(f_1)-(f_3)$. Assume that \eqref{c:intro} admits a unique positive solution, 
  then the following statements hold:
\begin{itemize}
  \item[$(1)$]  If $2+\frac{4}{N}<p<2^*$, then $\|u_\lambda\|_{L^2(\Omega)}^2\rightarrow 0$ as $\lambda \rightarrow +\infty$.
  \item[$(2)$] If $p=2+\frac{4}{N}$, then $\|u_\lambda\|_{L^2(\Omega)}^2\rightarrow \|Q\|_{L^2(\mathbb R^N)}^2$ as $\lambda \rightarrow +\infty$, 
  where $Q$ is the unique ground state of \eqref{global}.
  
  \item[$(3)$] If $2<p<2+\frac{4}{N}$, then $\|u_\lambda\|_{L^2(\Omega)}^2\rightarrow +\infty$ as $\lambda \rightarrow +\infty$.

\end{itemize}
  \end{theorem}

  \medskip
  
  If$\left ( \lambda, u_\lambda\right ) \in E$, then $\Phi:=e^{i\lambda t}u_\lambda$ is a standing wave solution of the  nonlinear 
  Schr\"{o}dinger equation (NLS)
  \begin{align}\label{orbit}
  \begin{cases}i\:\partial_t\Phi+\Delta \Phi+\tilde{f}(x, \Phi)=0,(t,x)\in(0,+\infty)\times\Omega,
    \\\Phi(0,x)=\Phi_{0}(x)=0,\end{cases}
  \end{align}
   where $\tilde{f}:\mathbb{C}\to\mathbb{C}$ is 
   defined by $\tilde{f}(x, e^{i\theta}u)=e^{i\theta}f(x, u),u\in\mathbb{R}.$ 
   Under the hypotheses of theorem~\ref{Normalizedsolution}, we can obtain 
   the orbital stability or instability of the standing waves associated with solutions of \eqref{orbit}. 
  
  \medskip
  
  We need some additional assumptions on \eqref{orbit}: 
  \begin{itemize}
    \item[$(LWP)$] For each $\Phi_0\in H_0^1(\Omega,\mathbb{C})$, there exists $t_0>0$
     only depending on $\|\Phi_0\|_{H_0^1(\Omega,\mathbb{C})}$, and a unique solution $\Phi(t,x)$ of ~\eqref{orbit}
    with initial datum $\Phi_0$ in the interval $I=[0,t_0)$.
    \item[$(N)$] Any positive solution of~\eqref{c:intro} is nondegenerate for all $\lambda \in (-\lambda_1,+\infty)$.

  \end{itemize}
  
  And we have the following result: 
  \begin{theorem}\label{orbitalstability}
    Let $N\geq 2$, $s\in (0,1)$, $\lambda> -\lambda_1(\Omega)$, $\Omega$ satisfy $(\Omega)$, $f$ satisfy $(f_1)-(f_3)$, and
    let $\Phi:=e^{i\lambda t}u_\lambda,\widetilde{\Phi}:=e^{i\lambda t}u_{\tilde{\lambda}}$, where $u_\lambda$ and
    $u_{\tilde{\lambda}}$ are given by Theorem~\ref{Normalizedsolution}. Assume that (LWP) and (N) holds. Then 
    the following statements hold:
  \begin{itemize}
    \item[$(1)$] 
    If $2+\frac{4}{N}<p<2^*$, for a.e. $c \in(0,b)$, $\Phi$ is orbitally stable while $\widetilde{\Phi}$ is orbitally unstable 
    where $-\lambda_1(\Omega)<\lambda<\tilde{\lambda}<+\infty$ and $\|u_\lambda\|_{L^2(\Omega)}^2=\|u_{\tilde{\lambda}}|\|_{L^2(\Omega)}^2=c$.
    
    \item[$(2)$] If $p=2+\frac{4}{N}$, then for a.e. $c \in(0,d)$, $\Phi$ is orbitally stable.

    \item[$(3)$] If $2<p<2+\frac{4}{N}$,
    for any $c \in(0,+\infty)$, $\Phi$ is orbitally stable, where $-\lambda_1(\Omega)<\lambda<+\infty$.
  
  \end{itemize}
  
    \end{theorem}

\begin{remark}
  When $s=1$, $f=|u|^{p-2}u$ and $p=2+\frac{4}{N}$, The authors in \cite{Noris1} proved
  that for any $0<c <\|Q\|_{L^2(\mathbb R^N)}^2$, there exists $(\lambda,u_\lambda) \in E$, such that $\|u_\lambda\|_{L^2(\Omega)}^2=c$; and
  any $c\geq\|Q\|_{L^2(\mathbb R^N)}^2$, there exists no $(\lambda,u_{\lambda})\in E$, such that $\|u_\lambda\|_{L^2(\Omega)}^2=c$; 
  and obtain the  orbital stability result based on this.
But when $f$ is mixed with $L^2$-critical and subcritical terms, or the equation is non-autonomous, 
  the more specific results may depend on more detailed analysis.
      
    \end{remark}
  
  We finish this section by introducing the structure of this paper.
  In section $2$, we will give some preliminary results.
  Section $3$ is devoted to the global bifurcation approach and the positive normalized solutions, that is, 
  proving Theorem \ref{Normalizedsolution} and \ref{uatinfity}.
  Then we apply our theorems to some examples in section $4$.
  Finally, we finish this paper via showing some orbital stability or instability results, i.e., Theorem \ref{orbitalstability}.

  \section{Preliminaries}

  Let $J_{\lambda}(u)$ be the energy functional given by

  $$
   J(u):=\frac12(\|\nabla u\|_{L^2(\Omega)}^2+\lambda\|u\|_{L^2(\Omega)}^2)-\int_\Omega F(x,u)dx, 
   $$
  
  then solutions of \eqref{c:intro} are critical points of $J_{\lambda}(u)$.
  
  For any $\lambda>-\lambda_1$ where $\lambda_1$ is the first eigenvalue of $-\Delta$ on $H_0^1(\Omega)$,
    we can define the Nehari manifold:
   
   $$
   \mathcal{N}:=\{u\in H_0^1(\Omega)\backslash\{0\}:\int_{\Omega}(|\nabla u|^2+\lambda u^2)dx-\int_{\Omega}f(x,u)udx=0\}.
   $$
   
  $$
  h(\lambda)=\inf_{N_{\lambda}}J_{\lambda}(u), 
  $$
  and
    $$
    K_{\lambda}=\{u\in N_{\lambda}: J_{\lambda}(u)=h(\lambda)\}.
    $$
  
    We know that $K_{\lambda}$ is nonempty for any $\lambda> -\lambda_1$ if $(f_1)$ and $(f_2)$ hold, and
    the elements in $K_{\lambda}$ are critical points of $J_{\lambda}$, which are called ground states.

We give some lemmas to describe some properties of the solutions of \eqref{c:intro}. 
\begin{lemma}
Assume that $(f_1)$ and $(f_2)$. Then the ground state of \eqref{c:intro} do not change sign.

\end{lemma}
\begin{proof}
By~\eqref{estimateA} we can verify that there exist $\phi>0$ and $\delta>0$ such that
 $\inf_{\int_{\Omega}|\nabla u|^2=\rho} J(v)>\delta>0$.
Thus for any $u \in H_0^1(\Omega) \backslash\{0\}$, there is a suitable $t$ with $t>0$ such that $J(tu)>\delta>0$. 
It follows from the proof of \cite[Lemma 4.1]{Minimax} that for any $u \in H_0^1(\Omega) \backslash\{0\}$, 
there exist a unique $t=t(u)$ such that $J(tu)$ achieves its maximum at $t(u)$.
As a consequence,  we deduce that 
$$
\inf_{u \in \mathcal{N}} J(u)>0.
$$
Let $u$ is a ground state of ~\eqref{c:intro}.
Assume by contradiction that $u$ is sign-changing. Define $u^+:=\max\{u,0\}$ and $u^-:=\min\{u,0\}$.
Then we know that $u^+,u^- \in \mathcal{N}$, $J(u^+)>0$, $J(u^-)>0$ and 
$J(u)=J(u^+)+J(u^-)$. 
On the other hand, since $u$ is a ground state, we have $J(u^+)\geq J(u)$ and $J(u^-)\geq J(u)$. 
Hence $J(u)=J(u^+)+J(u^-)\geq 2J(u)>0$, which yield a contradiction.

\end{proof}

\begin{lemma}\label{acurve}
If \eqref{c:intro} admits only one positive ground state $u_\lambda$ for any $\lambda>-\lambda_1$, then 
$K=\{(\lambda, u), \lambda>-\lambda_1, u \in K_\lambda, u>0\}$ is a continuous curve in $\mathbb R \times H_0^1(\Omega)$. 
Moreover, if $u_\lambda$ is the unique positive solution of ~\eqref{c:intro} for any $\lambda>-\lambda_1$,
and $u_\lambda \to 0 \text{ in } H_0^1(\Omega), \text{ as } \lambda \to -\lambda_1.$
  
  \end{lemma}
  
  \begin{proof}
Similar to the proof of~\cite[Corollary 2.4, Corollary 1.5]{LSongCVPDE}, we can obtain this result. 
  \end{proof}

  \begin{lemma}\label{usymmetry}
    Let $N\geq 2$, $s\in (0,1)$, $\lambda> -\lambda_1(\Omega)$. And let $f$ satisfy $(f_1)-(f_3)$, 
    $\Omega$ satisfy $(\Omega)$. If $u$ is a positive ground state of \eqref{c:intro}, then $u$ is $G$-invariant.
    
    \end{lemma}
  
\begin{proof}
  Since $f$ is symmetric in $x_i$ and nonincreasing in $|x_i|$ for all $i=1,\ldots,N$, 
  we can verify that $u$ is decreasing in $|x_i|$, symmetric in $x_i$ for all $i=1,\ldots,N$  by using the moving plane method.
  And then $u$ is $G$-invariant, by the rearrangement inequality. 
\end{proof}

\medskip

\section{A global bifurcation and normalized solutions}

Let $u_\lambda$ be a positive ground state, with respect to $\lambda$, of \eqref{c:intro}. 
We consider $u_\lambda$ when $\lambda > 1$ in the following arguments.

Let

$$
\mu=\frac{1}{\lambda}, 
$$
\begin{align}\label{definevmu}
v_\mu=\mu^{\frac{1}{p-2}}u_\lambda(\frac{x}{\lambda^{\frac{1}{2}}}),
 \end{align}
 and
$$
d_\mu(x)=d(\mu^\frac{1}{2}x).
$$

Then
$v_\mu>0$ on $\mu^{-\frac{1}{2}}\Omega=\{\mu^{-\frac{1}{2}}x : x \in \Omega\}$ (we denote $\mu^{-\frac{1}{2}}\Omega$ by $\Omega_\mu$), 
$v_\mu$ is $G$-invariant and when $\mu=1/\lambda\to0^+$, 
$\Omega_\mu \to \mathbb{R}^N$. 
Direct calculations yield that  
$$
\|v_\mu\|_{L^2(\Omega_\mu)}^2=\int_{\Omega_\mu}|\lambda^{-\frac{1}{p-2}}u_\lambda(\frac{x}{\lambda^{\frac{1}{2}}})|^2
=\int_{\Omega}\lambda^{\frac{N}{2}-\frac{2}{p-2}}|u_\lambda|^2
=\lambda^{\frac{N}{2}-\frac{2}{p-2}}\|u_\lambda\|_{L^2(\Omega)}^2,
$$

thus
$$
\|u_\lambda\|_{L^2(\Omega)}^2=\lambda^{\frac{2}{p-2}-\frac{N}{2}}\|v_\mu\|_{L^2(\Omega_\mu)}^2.
$$

Note that $v_\mu$ satisfies the following equation
 \begin{align}\label{newpeps}
     -\Delta v +v=\mu^{\frac{p-1}{p-2}}f(\mu^{\frac{1}{2}}x, \mu^{-\frac{1}{p-2}}v) 
     \qquad v=0\quad \text{ on }\partial \Omega_\mu.
     \end{align}

Choose the norm $\|v\|_\mu=\sqrt{\|\nabla v\|_{L^2(\Omega_\mu)}^2+\|v\|_{L^2(\Omega_\mu)}^2}$ for $H_0^1(\Omega_\mu)$.

We can verify that for any $\mu \in (0,1)$, $v_\mu$ is a positive ground state of equation \eqref{newpeps}, 
and the uniqueness and nondegeneracy of $v_\mu$ is equivalent to the the uniqueness and nondegeneracy of $u_\lambda$.

\begin{lemma}\label{lpandlploc}
  Assume that $\phi_n(x), \phi(x) \in L^\infty(\mathbb R^N)$,  
   $\phi_n(x) \to \phi(x)$ in $L^1_{loc}(\mathbb R^N)$, and $\psi(x) \in L^1(\mathbb R^N)$,  
  then 
  $$
  \int_{\mathbb R^N}(\phi_n(x)-\phi(x))\psi(x) \to 0,
  $$
  as $n \to +\infty$.
 
  \end{lemma}
  
  \begin{proof}
  Since $C_c^\infty(\mathbb R^N)$ is dense in $L^1(\mathbb R^N)$, then $\forall \epsilon>0$, there exists $\psi_0(x) \in C_c^\infty(\mathbb R^N)$
  such that 
  $$
  \int_{\mathbb R^N} |\psi_0(x)-\psi(x)|<\epsilon; 
  $$
  and $\exists N_0>0$, $\forall n>N_0$, 
  $$
  \int_{\mathbb R^N}|\phi_n(x)-\phi(x)\|\psi_0(x)|<\epsilon.
  $$
As a consequence, 
$$
 \begin{aligned}
  \int_{\mathbb R^N}(\phi_n(x)-\phi(x))\psi(x) &<\int_{\mathbb R^N}|\phi_n(x)-\phi(x)\|\psi_0(x)|+\int_{\mathbb R^N}|\phi_n(x)-\phi(x)|
  |\psi(x)-\psi_0(x)| \\
  &<(2C_1+1)\epsilon, 
  \end{aligned}
$$

  where $\|\phi_n(x)\|_{L^\infty(\mathbb R^N)}, \|\phi(x)\|_{L^\infty(\mathbb R^N)}<C_1.$

  \end{proof}

Let 
\begin{align}
  \widetilde{v}_\mu=
  \begin{cases}v_\mu, x \in \Omega_\mu,
      \\ 0, x \in \mathbb R^N\setminus \Omega_\mu, \end{cases}
  \end{align} 
  be the zero extension function in $\mathbb R^N$ of $v_\mu \in H_0^1(\Omega_\mu)$, then $\widetilde{v}_\mu \in H_G^1(\mathbb R^N)$.
  Similarly, let $\widetilde{d}_\mu$  be the zero extension function of $d_\mu$ in $\mathbb R^N$, 
then $\widetilde{d}_\mu \to d(0)=1$ as $\mu \to 0$ in $L^1_{loc}(\mathbb R^N)$.

Define
$$
 \begin{aligned}\label{definephi}
  \Phi_{\mu}(v) &:=\frac{1}{2}(\|\nabla v\|_{L^2(\Omega_\mu)}^2+\|v\|_{L^2(\Omega_\mu)}^2)
   -\mu^{\frac{p}{p-2}}\int_{\Omega_\mu}F(\mu^{\frac{1}{2}}x, \mu^{-\frac{1}{p-2}}v) \\
  &=\frac{1}{2}(\|\nabla v\|_{L^2(\Omega_\mu)}^2+\|v\|_{L^2(\Omega_\mu)}^2)
   -\frac{1}{p}\int_{\Omega_\mu}d_\mu(x)|v|^p-\mu^{\frac{p}{p-2}}\int_{\Omega_\mu}G(\mu^{\frac{1}{2}}x, \mu^{-\frac{1}{p-2}}v), 
  \end{aligned}
$$
$$
\widetilde{\mathcal{N}}_{\mu}:=\{v \in H_0^1(\Omega_\mu)\setminus {0} :\: \Phi'_\mu(v)v=0\},
$$
$$
\widetilde{h}(\mu):=\inf_{\widetilde{\mathcal{N}}_{\mu}}\Phi_{\mu}(v), 
$$

for $0<\mu<1$; 

and define
 $$
 \Phi_{0}(v):=\frac12(\|\nabla v\|_{L^2(\mathbb R^N)}^2+\|v\|_{L^2(\mathbb R^N)}^2)
 -\frac{1}{p}\int_{\mathbb R^N}|v|^p,
  $$
  $$
  \widetilde{\mathcal{N}}_{0}:=\{v \in H^1(\mathbb R^N)\setminus {0} :\: \Phi'_0(v)v=0\},
  $$
  and
  $$
\widetilde{h}(0):=\inf_{\widetilde{\mathcal{N}}_{0}}\Phi_{0}(v).
$$

 Note that $\Phi_0(v)$ is the energy functional of the following equation:
$$
  -\Delta v +v=|v|^{p-2}v \quad\mbox{ in } \mathbb R^N.
$$
     
  Choose the norm $\|v\|=\sqrt{\|\nabla v\|_{L^2(\mathbb R^N)}^2+\|v\|_{L^2(\mathbb R^N)}^2}$ for $H^1(\mathbb R^N)$.

  \medskip

  Define $a:= \frac{1}{2}dist(0, \partial\Omega)$, and $M_1:=\{x \in \Omega :x+a\frac{x}{|x|} \in \Omega\}$. 
  Here we do not distinguish between $x$ as a vector or as a point, 
  since they belong to an equivalence class and do not lead to confusion.
  Let $\gamma_1 \in C_c^\infty(\mathbb R)$ be a truncation function such that  $\gamma_1(x)=1$ for $x\in M_1$, 
  $0< \gamma_1\leq 1$ for $x \in \Omega \backslash M_1$, and
 $\gamma_1(t)=0$ in $\mathbb R^N \backslash \Omega$.
  Define $M_\mu:=\{x \in \Omega_\mu :x+a\frac{x}{|x|} \in \Omega_\mu\}$.
  Let $\gamma_\mu \in C_c^\infty(\mathbb R^N)$  such that $\gamma_\mu(x)=1$ for $x \in M_\mu$, 
  $\gamma_\mu(x)=\gamma_1(x-(1-\mu^{\frac{1}{2}})x_0)$ for $x \in \Omega_\mu \backslash M_\mu$,
  where $x_0 \in \partial \Omega_\mu$ and $x_0=kx$ for fixed $x \in \Omega_\mu \backslash M_\mu$ and $k>0$, 
  (indeed, the shape of domain ensures the existence and uniqueness of $x_0$,) 
  and $\gamma_\mu(x)=0$ for $x \in \mathbb R^N \backslash \Omega_\mu$.

Recall that \eqref{global} admits a unique positive ground state $Q$, and $Q \in H_G^1(\mathbb R^N)$. 
Let $\widetilde{Q}_\mu=\gamma_\mu Q$, and 
$Q_\mu=(\widetilde{Q}_\mu)_{|\Omega_\mu}$,  
by the definition of $\gamma_\mu$ we know that $\widetilde{Q}_\mu \in H_G^1(\mathbb R^N)$, and
$Q_\mu \in H_{0,G}^1(\Omega_\mu)$. 
 It is not difficult to verify that $Q_\mu=\widetilde{Q}_\mu$ as $\mu \to 0^+$ and
  $\widetilde{Q}_\mu \to Q$ in $H_G^1(R^N)$ as $\mu \to 0^+$, thus we can define $Q_0=Q$ in the following arguments.

 \begin{lemma}\label{hmubound}
  Let $N\geq 2$, $\lambda> -\lambda_1$, $\Omega$ satisfy $(\Omega)$ and $f$ satisfy $(f_1)-(f_3)$. Then
  there exists $\delta=\delta(Q)>0$ and a constant $C_1=C_1(Q)>0$ such that $h(\mu)<C_1$ for $0<\mu<\delta$.
  
    \end{lemma}

    \begin{proof}

   First we prove that, for any $\mu \in [0,1)$, there exists a unique $t:[0,1)  \to (0,+\infty)$ 
    such that $tQ_\mu \in \widetilde{\mathcal{N}}, t>0 \Leftrightarrow t=t(\mu)$, where
    $$
    \widetilde{\mathcal{N}}=\cup_{\mu \in [0,1)}\widetilde{\mathcal{N}}_\mu.
    $$

    For $0<\mu<1$, Consider
     $$
     \begin{aligned}
      \psi: & (0,\infty)\times(0,1) \rightarrow \mathbb{R},\\
      & (t,\mu)\mapsto \Phi_\mu(tQ_\mu),
    \end{aligned}
    $$
  i.e.
$$
\psi(t,\mu)=\frac{1}{2} t^2\|\nabla Q_\mu\|_{L^2(\Omega_\mu)}^2+\frac{1}{2}t^2\|Q_\mu\|_{L^2(\Omega_\mu)}^2-
\mu^{\frac{p}{p-2}}\int_{\Omega_\mu} F(\mu^{\frac{1}{2}}x, (\mu^{-\frac{1}{p-2}}tQ_\mu)) dx.
$$

Then
\begin{align}
 \psi_{t}(t,\mu)=t\|\nabla Q_\mu\|_{L^2(\Omega_\mu)}^2+t\|Q_\mu\|_{L^2(\Omega_\mu)}^2-
 \mu^{\frac{p-1}{p-2}}\int_{\Omega_\mu} f(\mu^{\frac{1}{2}}x, (\mu^{-\frac{1}{p-2}}tQ_\mu))Q_\mu dx.
\end{align}
     
Note that $(\mu,tQ_\mu)\in \widetilde{\mathcal{N}}$ if and only if $\psi_t(t,\mu)=0$,
i.e.
 \begin{align}\label{uniquetu}
 \|\nabla Q_\mu\|_{L^2(\Omega_\mu)}^2+\|Q_\mu\|_{L^2(\Omega_\mu)}^2=
 \frac{1}{t}\mu^{\frac{p-1}{p-2}}\int_{\Omega_\mu} f(\mu^{\frac{1}{2}}x, (\mu^{-\frac{1}{p-2}}tQ_\mu))Q_\mu dx. 
 \end{align}
    
When $\psi_t(t,\mu)=0$, we have

  $$
  \begin{aligned}
    \psi_{tt}(t,\mu) &=\frac{1}{t}\mu^{\frac{p-1}{p-2}}\int_{\Omega_\mu} f_u(\mu^{\frac{1}{2}}x,(\mu^{-\frac{1}{p-2}}tQ_\mu))Q_\mu dx
    -\int_{\mathbb R^N} f_u(\mu^{\frac{1}{2}}x,(\mu^{-\frac{1}{p-2}}tQ_\mu))(Q_\mu)^2 dx \\
    &=\frac{\int_{\Omega_\mu} \mu^{\frac{p-1}{p-2}}Q_\mu [f_u(\mu^{\frac{1}{2}}x,(\mu^{-\frac{1}{p-2}}tQ_\mu))-
    f_u(\mu^{\frac{1}{2}}x,(\mu^{-\frac{1}{p-2}}tQ_\mu)(t\mu^{-\frac{1}{p-2}}Q_\mu))]}{t}\\
    &<0.
    \end{aligned}
$$

By \eqref{estimateB} and~\eqref{Q:intro3B}, there exists some $r>0$ small enough such that
\begin{align}\label{Phiuniformge0}
  \inf_{||\widetilde{v}\|=r}\Phi_{\mu}(v)>\frac{1}{2}r^2-Cr^p-O(r^2)>0, 
\end{align}

where $O(r^2)$ means that $\lim_{r \to 0}\frac{O(r^2)}{r^2}=0$. 
Then we know that for some suitable $t$, $\psi(t,\mu)\geq\frac{1}{2}r^2-Cr^p-O(r^2)>0$. 
By \eqref{estimateB}, we derive that $\psi(t,\mu) \to -\infty$ as $t \to +\infty$ for fixed $\mu \in (0,1)$.
Thus we get the existence and uniqueness of $t(\mu)$.
  
For $\mu=0$, the proof is similar, and we can deduce that $t(0)=1$. We will show that $t(\mu) \to t(0)$ as $\mu \to 0^+$. 
  
Note that
\begin{eqnarray} \label{split}
&&  
\frac{1}{t}\mu^{\frac{p-1}{p-2}}\int_{\Omega_\mu} f_u(\mu^{\frac{1}{2}}x, (\mu^{-\frac{1}{p-2}}tQ_\mu))Q_\mu dx=
\nonumber\\
&&  
=t^{p-1}(\int_{\Omega_\mu}d_\mu(x)|Q_\mu|^p+
  \frac{1}{t^p}\mu^{\frac{p-1}{p-2}}\int_{\Omega_\mu} g_u(\mu^{\frac{1}{2}}x, (\mu^{-\frac{1}{p-2}}tQ_\mu))Q_\mu dx).
\nonumber\\
\end{eqnarray}

Together with the uniqueness of $t(\mu)$, observe \eqref{uniquetu} and \eqref{split} again we see that,  
when $\mu \to 0^+$, 
$\|\nabla Q_\mu\|_{L^2(\Omega_\mu)}^2$ and $\|Q_\mu\|_{L^2(\Omega_\mu)}^2$  converge to some constants not equaling to $0$.
By Lemma~\ref{lpandlploc}, we have
\begin{align}\label{dandQmuto0}
\int_{\Omega_\mu}d_\mu(x)|Q_\mu|^p=\int_{\mathbb R^N} \widetilde{d}_\mu(x) |\widetilde{Q}_\mu|^p \to 
\int_{\mathbb R^N} |Q(x)|^p.
\end{align}

Indeed, since $\widetilde{d}_\mu(x) \to 1$ in $L_{loc}^{1}(\mathbb R^N)$, 
$$
\int_{\mathbb R^N} \widetilde{d}_\mu(x) |\widetilde{Q}_\mu|^p-|Q(x)|^p
= \int_{\mathbb R^N}(\widetilde{d}_\mu(x)-1)|Q|^p
+\int_{\mathbb R^N}\widetilde{d}_\mu(x)(|\widetilde{Q}_\mu|^p-|Q(x)|^p) \to 0.
$$

Assume that there exists $\{\mu_n\} \subset (0,1)$ with $\mu_n \to 0$ as $n \to +\infty$, 
such that $t(\mu_n) \to +\infty$ as $n \to +\infty$. By \eqref{Q:intro3B}, we have
\begin{eqnarray}\label{gtermto0}
  &&  
  \frac{1}{(t(\mu_n))^p}\mu_n^{\frac{p-1}{p-2}}\int_{\Omega_{\mu_n}} g_u({\mu_n}^{\frac{1}{2}}x, 
  ({\mu_n}^{-\frac{1}{p-2}}tQ_{\mu_n}))Q_{\mu_n} dx=
  \nonumber\\
  &&  
  =\frac{1}{(t(\mu_n))^p}\frac{1}{({\mu_n}^{-\frac{1}{p-2}})^{p-1}}\int_{\Omega_{\mu_n}} g_u({\mu_n}^{\frac{1}{2}}x, ({\mu_n}^{-\frac{1}{p-2}}
  tQ_{\mu_n}))Q_{\mu_n} dx \to 0.
  \nonumber\\
  \end{eqnarray}
Combining~\eqref{split}, ~\eqref{dandQmuto0} and~\eqref{gtermto0}, we obtain that
$$
\frac{1}{t(\mu_n)}\mu^{\frac{p-1}{p-2}}\int_{\Omega_{\mu_n}} f_u(|{\mu_n}^{\frac{1}{2}}x|, ({\mu_n}^{-\frac{1}{p-2}}
tQ_{\mu_n}))Q_{\mu_n} dx \to +\infty
$$
as $n \to \infty$. Contradiction!
  
Therefore, there is no sequence $\{\mu_n\} \subset (0,1)$ with $\mu_n \to 0$ as $n \to \infty$ 
such that $\{t(\mu_n)\}$ is unbounded, i.e. $t(\mu)$ is  uniformly bounded with respect to $\mu$.
Thus we know that $\widetilde{h}(\mu)$ is bounded as $\mu \to 0^+$ since $\widetilde{h}(\mu) \leq \Phi_\mu(t(\mu)Q_\mu)$. 

Since $t(\mu)$ is bounded as $\mu \to 0^+$, and $Q_\mu \to Q$ in $H^1(\mathbb R^N)$, 
if a subsequence $\{\mu_n\} \subset (0,1)$ converging to $0$ such that $t(\mu_n) \to t_0$,
it follows from \eqref{uniquetu} that $t_0=t(0)=1$. Hence $t(\mu) \to 1$ as $\mu \to 0^+$.
Letting $C_1=\Phi_0(Q)+1$, we can take $\delta>0$ small enough such that $h(\mu)<C_1$ for $0<\mu<\delta$.
\end{proof}

  \begin{proposition}\label{vmubound}
   Under the hypotheses of lemma~\ref{hmubound}. Assume $v_\mu$ is the defined in \eqref{definevmu}, as $\mu \to 0^+$, then
   there exists a positive constant $C=C(Q)$, $\|v_\mu\|_\mu \leq C$.

    \end{proposition}

    \begin{proof}
    By lemma \ref{hmubound}, there exists $\delta=\delta(Q)>0$ such that $h(\mu)<C_1$ for $0<\mu<\delta$.
    Together with $v_\mu \in \widetilde{K}_\mu$, we have
    
    $$
    \begin{cases}\Phi_{\mu}(v_\mu):=\frac{1}{2}(\|\nabla v_\mu\|_{L^2(\Omega_\mu)}^2+\|v_\mu\|_{L^2(\Omega_\mu)}^2)
      -\mu^{\frac{p}{p-2}}\int_{\Omega_\mu}F(\mu^{\frac{1}{2}}x, \mu^{-\frac{1}{p-2}}v_\mu)<C_1,
        \\\langle\Phi_{\mu}^{\prime}(v_\mu),v_\mu\rangle=\|\nabla v_\mu\|_{L^2(\Omega_\mu)}^2+\|v_\mu\|_{L^2(\Omega_\mu)}^2-
        \mu^{\frac{p-1}{p-2}}\int_{\Omega_\mu} f(\mu^{\frac{1}{2}}x, \mu^{-\frac{1}{p-2}}v_\mu)v_\mu=0.\end{cases}
    $$

  Then we have
  $$
  \frac{1}{2}\mu^{\frac{p-1}{p-2}}\int_{\Omega_\mu} f(\mu^{\frac{1}{2}}x, \mu^{-\frac{1}{p-2}}v_\mu)v_\mu - 
  \mu^{\frac{p}{p-2}}\int_{\Omega_\mu}F(\mu^{\frac{1}{2}}x, \mu^{-\frac{1}{p-2}}v_\mu) <C_1,
  $$
    
   Since $0< \alpha F(\mu^{\frac{1}{2}}x, \mu^{-\frac{1}{p-2}}v_\mu)
   \leq f(\mu^{\frac{1}{2}}x, \mu^{-\frac{1}{p-2}}v_\mu)\mu^{-\frac{1}{p-2}}v_\mu$,
  it follow that 
  $$
  \mu^{-\frac{p-1}{p-2}}\int_{R^N}f(\mu^{\frac{1}{2}}x, \mu^{-\frac{1}{p-2}}v_\mu)v_\mu<\frac{2\alpha}{\alpha-2}C_1,
  $$
  thus $\|v_\mu\|_\mu=\sqrt{\|\nabla v_\mu\|_{L^2(\Omega_\mu)}^2+\|v_\mu\|_{L^2(\Omega_\mu)}^2}<\sqrt{\frac{2\alpha}{\alpha-2}C_1}.$

    \end{proof}

\begin{proposition}\label{vandQ}
  Under the hypotheses of lemma~\ref{hmubound}, there holds that
 $\widetilde{v}_\mu \to Q$ in $H_G^1(\mathbb R^N)$ as $\mu \to 0^+$, where $Q$ is the unique positive solution of \eqref{global}.
\end{proposition}

\begin{proof}

Consider $t(\mu)Q_\mu \in \widetilde{\mathcal{N}}$, since $t(\mu)\widetilde{Q}_\mu\to Q$ in $H_G^1(\mathbb R^N)$, 
we have

$$
  \begin{aligned}
    \widetilde{h}(\mu) &\leq \Phi_\mu(t(\mu)Q_\mu) \\
    &=\Phi_0(t(\mu)\widetilde{Q}_\mu)+\frac{1}{p}\int_{\mathbb R^N}(t(\mu)\widetilde{Q}_\mu)^{p}-\mu^{\frac{p}{p-2}}
    \int_{\Omega_\mu}F(\mu^{\frac{1}{2}}x, (\mu^{-\frac{1}{p-2}}t(\mu)Q_\mu))\\
    &=\Phi_0(Q)+o(1)\\
    &=\widetilde{h}(0)+o(1).
    \end{aligned}
$$

Letting $\mu \to 0^+$, we have 
\begin{align}\label{hmuleqhh0}
\limsup_{\mu \to 0^+}\widetilde{h}(\mu) \leq \widetilde{h}(0).
\end{align}

Take an arbitrary sequence $\{\mu_n\} \subset (0,1)$ such that $\mu_n \to 0$ as $n \to \infty$. 
By Proposition~\ref{vmubound}, we can assume $\widetilde{v}_{\mu_n}$ is uniformly bounded in $H^1(\mathbb R^N)$ without loss of generality . 
Let $v_n=v_{\mu_n}$, $\widetilde{v}_n=\widetilde{v}_{\mu_n}$, $d_n(x)=d_{\mu_n}(x)$ and $\widetilde{d}_n(x)=\widetilde{d}_{\mu_n}(x)$. 

Direct calculations yield that 
\begin{align}\label{tGto0}
  \mu_n^{\frac{p}{p-2}}\int_{\Omega_{\mu_n}}G(\mu_n^{\frac{1}{2}}x, \mu_n^{-\frac{1}{p-2}}v_n) \to 0,
\end{align}
as $n \to \infty$. 

Moreover, it follows from Lemma~\ref{dn-1vnto0} that
\begin{align}\label{dvto0}
  \frac{1}{p}\int_{\Omega_{\mu_n}}d_n(x)|v_n|^p-\frac{1}{p}\int_{\Omega_{\mu_n}}|v_n|^p=
\frac{1}{p}\int_{\mathbb R^N}(\widetilde{d}_n(x)-1)|\widetilde{v}_n|^p \to 0, 
\end{align}
as $n \to \infty$. 

Combining~\eqref{tGto0} and~\eqref{dvto0}, we have
\begin{align}\label{F-vnpto0}
  {\mu_n}^{\frac{p}{p-2}}\int_{\Omega_{\mu_n}}F(\mu_n^{\frac{1}{2}}x, \mu_n^{-\frac{1}{p-2}}v_n)-
\frac{1}{p}\int_{\Omega_{\mu_n}}|v_n|^p \to 0,
\end{align}

and
$$
{\mu_n}^{\frac{p}{p-2}}\int_{\Omega_{\mu_n}}F({\mu_n}^{\frac{1}{2}}x, {\mu_n}^{-\frac{1}{p-2}}v_n)-
\frac{1}{p}\int_{\mathbb R^N}|\widetilde{v}_n|^p \to 0 ,
$$ 
as $n \to \infty$.

Therefore, 
\begin{align}\label{phi0vn}
\Phi_{0}(\widetilde{v}_n)=\Phi_{{\mu_n}}(v_n)+{\mu_n}^{\frac{p}{p-2}}\int_{\Omega_{\mu_n}}F(\mu_n^{\frac{1}{2}}x, \mu_n^{-\frac{1}{p-2}}v_n)
-\frac{1}{p}\int_{\mathbb R^N}{\widetilde{v}_n}^p
=\widetilde{h}({\mu_n})+o(1).
\end{align}

Together with~\eqref{hmuleqhh0}, we can derive that 
\begin{align}\label{phi0leh0}
\limsup_{n \to \infty} \Phi_{0}(\widetilde{v}_n) \leq \widetilde{h}(0).
\end{align}

By~\eqref{Phiuniformge0}, we have
$$
\liminf_{n \to \infty}\widetilde{h}(\mu_n)\geq \liminf_{n \to \infty} \inf_{\|v\|=r}\Phi_{\mu_n}(v)>\frac{1}{2}r^2-Cr^p-O(r^2)>0.
$$

Therefore, by~\eqref{phi0vn} and~\eqref{phi0leh0}, we can assume that $\Phi_{0}(\widetilde{v}_n) \to c \in (0, \widetilde{h}(0)]$ 
up to a subsequence if necessary.

Similar to the construction of $Q_\mu$, for any $\phi \in H_G^1(\mathbb R^N)$, we can find
$\widetilde{\phi}_{\mu_n} \in H_G^1(\mathbb R^N)$ and
$\phi_{\mu_n} \in H_{0,G}^1(\Omega_{\mu_n})$
(denoted by $\widetilde{\phi}_n$ and $\phi_{n}$) such that $\widetilde{\phi}_n \to \phi$ in $H_G^1(\mathbb R^N)$.
For any $\phi \in H_G^1(\mathbb R^N)$, we have
$$
\begin{aligned}
  \langle\Phi_{0}'(\widetilde{v}_n), \phi\rangle &=\int_{\mathbb R^N}\nabla \widetilde{v}_n \nabla \phi 
  +\int_{\mathbb R^N} \widetilde{v}_n \phi-\int_{\mathbb R^N}|\widetilde{v}_n|^{p-2}\widetilde{v}_n \phi\\
  &=\int_{\mathbb R^N}\nabla \widetilde{v}_n \nabla \widetilde{\phi}_n 
  +\int_{\mathbb R^N} \widetilde{v}_n \widetilde{\phi}_n
  -\int_{\mathbb R^N}|\widetilde{v}_n|^{p-2}\widetilde{v}_n \phi\\
  &\qquad+\int_{\mathbb R^N}\nabla \widetilde{v}_n \nabla (\phi-\widetilde{\phi}_n) 
  +\int_{\mathbb R^N} \widetilde{v}_n (\phi-\widetilde{\phi}_n),
  \end{aligned}
$$

and
$$
\langle\Phi_{\mu_n}'(v_n), \phi_n\rangle
=\int_{\Omega_{\mu_n}}\nabla v_n \nabla \phi_n +\int_{\Omega_{\mu_n}} v_n \phi_n
-\int_{\Omega_{\mu_n}}{\mu_n}^{\frac{p-1}{p-2}}f(\mu_n^\frac{1}{2}x, {\mu_n}^{-\frac{1}{p-2}}v_n) \phi_n. 
$$

Then by \eqref{F-vnpto0}, we derive that 
$$
  \begin{aligned}
    \langle\Phi_{0}'(\widetilde{v}_n), \phi\rangle &=\langle\Phi_{\mu_n}'(v_n), \phi_n\rangle
    +\langle\Phi_{0}'(\widetilde{v}_n), \phi-\widetilde{\phi}_n\rangle
    +\int_{\mathbb R^N}|\widetilde{v}_n|^{p-2}\widetilde{v}_n (\phi-\widetilde{\phi}_n)\\
    &\qquad+\int_{\Omega_{\mu_n}}{\mu_n}^{\frac{p-1}{p-2}}f(\mu_n^\frac{1}{2}x, {\mu_n}^{-\frac{1}{p-2}}v_n) \phi_n
    -\int_{\mathbb R^N}|\widetilde{v}_n|^{p-2}\widetilde{v}_n \phi\\
  &=\langle\Phi_{\mu_n}'(v_n), \phi_n\rangle+\langle\Phi_{0}'(\widetilde{v}_n), \phi-\widetilde{\phi}_n\rangle\\
  &\qquad+\int_{\Omega_{\mu_n}}\left({\mu_n}^{\frac{p}{p-2}}f(\mu_n^\frac{1}{2}x, {\mu_n}^{-\frac{1}{p-2}}v_n)
  -|v_n|^{p-2}v_n\right)\phi_n\\
  &=\langle\Phi_{\mu_n}'(v_n), \phi_n\rangle+\langle\Phi_{0}'(\widetilde{v}_n), \phi-\widetilde{\phi}_n\rangle
  +\langle H'(\mu_n, v_n),\phi_n \rangle\\
  & \to 0
    \end{aligned}
$$
as $n \to \infty$, where 
$$
H(\mu_n, v_n)=\int_{\Omega_{\mu_n}}({\mu_n}^{\frac{p}{p-2}}F(\mu_n^{\frac{1}{2}}x, \mu_n^{-\frac{1}{p-2}}v_n)-\frac{1}{p}|v_n|^p ).
$$
Thus ${\Phi_0}_{|H_G^1(\mathbb R^N)}'(\widetilde{v}_n) \to 0$ as $n \to +\infty$.

Since $\widetilde{v}_n \in H_G^1(\mathbb R^N)$, and embedding $H_G^1(\mathbb R^N) \hookrightarrow L^{q}(\mathbb{R}^N)$ is compact 
for $q\in(2,2^*)$, 
we know that ${\Phi_0}_{|H_G^1(\mathbb R^N)}$ satisfies $(PS)_c$ condition for $c \in (0, \widetilde{h}(0)]$, 
then $\widetilde{v}_n \to v_0$ in $H_G^1(\mathbb R^N)$, and ${\Phi_0}'_{|H_G^1(\mathbb R^N)}(v_0)=0$.

Note that
$$
\Phi_{0}(v_0)=\lim_{n \to \infty} \Phi_{0}(\widetilde{v}_n)>0,
$$
we know that $v_0 \neq 0$. Hence $v_0 \geq 0$ in $\mathbb R^N$ and $v_0$ is a nontrivial solution of~\eqref{global}. 
By the strong maximum principle, $v_0>0$ and $v_0$ is a positive solution of~\eqref{global}. 
It follows from the uniqueness of the positive solution of \eqref{global} that $v_0=Q$.

\end{proof}

\begin{remark}
If $Q$ is the unique positive ground state of~\eqref{global}, we can still complete the proof of Proposition~\ref{vandQ}
even we do not know whether the positive solution of ~\eqref{global} is unique. 
Indeed, we can prove that $v_0$ is a ground state of \eqref{global}.

By \eqref{phi0vn}, we deduce that $\forall \epsilon >0, \exists N_0>0, \forall n>N_0$, there holds that

  $$
  \begin{aligned}
    \widetilde{h}(0) &\leq \Phi_{0}(v_0) \\
    &\leq \Phi_{0}(\widetilde{v}_n)+\frac{\epsilon}{2}\\
    &=\Phi_{\mu_n}(\widetilde{v}_n)+\frac{1}{p}\int_{\mathbb R^N}{|\widetilde{v}_n|}^p
    -{\mu_n}^{\frac{p-1}{p-2}}\int_{\Omega_\mu}F({\mu_n}^{-\frac{1}{p-2}}v_n)+\frac{\epsilon}{2}\\
    &\leq \Phi_{{\mu_n}}(\widetilde{v}_n) + \epsilon.
    \end{aligned}
$$

Letting $n \to \infty$, we know that 
$$
\widetilde{h}(0) \leq \limsup_{n \to \infty}\widetilde{h}({\mu_n}).
$$
Therefore,  
$$
\widetilde{h}(0) = \limsup_{\mu\to 0^+}\widetilde{h}({\mu}).
$$
We derive that $\Phi_{0}(v_0)=\widetilde{h}(0)$, thus $v_0$ is a positive ground state of \eqref{global}.
By the uniqueness of the positive ground state of \eqref{global}, $v_0=Q$.

\end{remark}

\begin{lemma}\label{dn-1vnto0}
Assume $\widetilde{d}_n(x)$ and $\widetilde{v}_n$ is defined as above, then
\begin{align}\label{formualdn-1nvto0}
  \int_{\mathbb R^N}(\widetilde{d}_n(x)-1)\widetilde{v}_n^p \to 0.
\end{align}

\end{lemma}

\begin{proof}
Since $\widetilde{v}_n$ is uniformly bounded in $H^1(\mathbb R^N)$, 
then by embedding $H^1(\mathbb R^N) \hookrightarrow L^p(\mathbb R^N)$, 
$\widetilde{v}_n$ is uniformly bounded in $L^p(\mathbb R^N)$, where $2\leq p<2^*$.
Therefore, $\widetilde{v}_n^p$ is uniformly bounded in $L^1(\mathbb R^N)$.
By Lemma~\ref{usymmetry},  $\widetilde{v}_n^p$ is decreasing in $|x_i|$ for all $i=1,\ldots,N$. 

For any $\epsilon>0$, we claim that there is some $k>0$ such that $D=k \Omega$ satisfies 
\begin{align}\label{uniformforvnp}
  \int_{\mathbb R^N \backslash D}\widetilde{v}_n^p <\epsilon~~~\text{uniformly}. 
\end{align}

Assume by contradiction that for any $D_n=n\Omega$ there is a subsequence of $\widetilde{v}_n$ (still denoted by $\widetilde{v}_n$, )
such that $\int_{D_n}\widetilde{v}_n^p \geq \epsilon$. 
Since $\widetilde{v}_n^p$ is decreasing in $|x_i|$ for all $i=1,\ldots,N$, there must be $\delta=\delta(\epsilon)$ 
such that $\widetilde{v}_n^p \geq \delta$ on $\partial D_n$ for any $n$. 
Then 
$$
\int_{\mathbb R^N}\widetilde{v}_n^p \geq \int_{D_n}\widetilde{v}_n^p \geq \delta|D_n| \to +\infty~~\text{as}~~n \to \infty, 
$$
yielding a contradiction.

Then we have
$$
\int_{\mathbb R^N \backslash D}|\widetilde{d}_n(x)-1|\widetilde{v}_n^p \leq \int_{\mathbb R^N \backslash D}\widetilde{v}_n^p<\epsilon. 
$$
On the other hand, for the same $\epsilon$ in~\eqref{uniformforvnp}, 
there exists $X=X(\epsilon)$ such that  $|\widetilde{d}_n(x)-1|<\epsilon$ for $n>X$ on $D$. Thus
$$
\int_{D} |\widetilde{d}_n(x)-1|\widetilde{v}_n^p < C_0 \epsilon, 
$$
where $C_0$ is the uniform $L^1(\mathbb R^N)$ bound of $\widetilde{v}_n^p$.

Note that 
$$
\int_{\mathbb R^N}|\widetilde{d}_n(x)-1|\widetilde{v}_n^p 
=\int_{D} |\widetilde{d}_n(x)-1|\widetilde{v}_n^p + \int_{\mathbb R^N \backslash D}|\widetilde{d}_n(x)-1|\widetilde{v}_n^p,
$$
we have $\int_{\mathbb R^N}|\widetilde{d}_n(x)-1|\widetilde{v}_n^p \leq (C_0+1) \epsilon$. 
From the arbitrariness of $\epsilon$, we obtain~\eqref{formualdn-1nvto0}.

\end{proof}

\begin{proof} [\textbf{Proof of Theorem~\ref{uatinfity}}]

  Since $\|u_\lambda\|^2_{L^2(\Omega)}=\lambda^{\frac{2}{p-2}-\frac{N}{2}}\|v_\mu\|^2_{L^2(\Omega_\mu)}$, 
 by proposition~\ref{vandQ}, we complete the proof.

\end{proof}

Then we can prove theorem \ref{Normalizedsolution}.
\begin{proof}[\textbf{Proof of Theorem~\ref{Normalizedsolution}}]

According to Lemma~\ref{acurve}, $\|u_\lambda\|^2_{L^2(\Omega)}$ is continuous 
with $\lambda$ in $(-\lambda_1,+\infty)$, and 
$u_\lambda \to 0 \text{ in } H_0^1(\Omega), \text{ as } \lambda \to -\lambda_1.$

By Theorem~\ref{Normalizedsolution}, for $2+\frac{4}{N}<p<2^*$, 
let $b=\max_{\lambda \in (-\lambda_1,+\infty)}\|u_\lambda\|^2_{L^2(\Omega)}$; 
and for $p=2+\frac{4}{N}$, 
let $D=\max_{\lambda \in (-\lambda_1,+\infty)}\|u_\lambda\|^2_{L^2(\Omega)}$,
thus $D \geq d=\|Q\|_{L^2(\mathbb R^N)}^2>0$. 
Therefore, we complete the proof of Theorem~\ref{Normalizedsolution}.
\end{proof}

\section{Some examples}

Consider the equation 
\begin{align}\label{equationhr}
  \begin{cases}-\Delta u +\lambda u= h(r) u^{p-1} \quad \text{ in }B_R,
    \\ u>0, \quad u_{|\partial B}=0, \end{cases}
\end{align}
where $N\geq 2$, $B_R$ is a ball with radius $R$, and $p \in (2,2^*)$. 
We introduce a theorem from \cite[Theorem 2.1]{Yanagida}.

\begin{theorem}\label{uniquenessfromE}

Let $m\in [0,n-2]$ be a parameter and define 
$$
    \begin{aligned}
    H(r;m)& :=2r^{m+2} h_{r}(r)/(p) \\
    &-\left\{2n-4-m-2(m+2)/(p)\right\}r^{m+1}h(r).
    \end{aligned}
$$

Assume that the following conditions hold:
  \begin{itemize}
  
    \item[$(C_1)$] $h(r)\geq0$ for all $r\in(0,R)$ and $h(r)>0$ for some $r\in(0,R)$, where $R$ is the radius of $B$.
  
    \item[$(C_2)$] $H(r;0)=2r^2h_r(r)/(p)-2\{n-2-2/(p)\}rh(r)\leq0$ for all
    $r\in(0,R).$
  
    \item[$(C_3)$] For each $m\in(0,n-2]$, there exists a $\beta(m)\in[0,R]$ such that $H(r;m)\geq0$
    for $r\in(0,\beta(m))$ and $H(r;m)\leq0$ for $r\in(\beta(m),R).$

   Then \eqref{equationhr} admits at most one solution. 

  \end{itemize}

\end{theorem}

We give an additional assumption on $h$:
\begin{itemize}
  \item[$(C_4)$] 
$h(0)>0$, and there exists $M_0$ such that $h(r)< M_0$ for any $r\in [0,R)$. 
\end{itemize}
Indeed, we can obtain the existence of solutions of \eqref{equationhr} by standard variational methods if $(C_1)$ and $(C_4)$ hold. 
Without loss of generality, we assume that $h(0)=1$. 
Therefore, we have the following theorem:

\begin{theorem}
  Let $N\geq 2$, $\lambda> -\lambda_1$. Assume that $(C_1)$-$(C_4)$ hold, then the following statements hold.
\begin{itemize}
\item[$(1)$] 
If $2+\frac{4}{N}<p<2^*$, then there exists some $b>0$ such that 
\begin{itemize}
\item[$(i)$] there exists at least one normalized solution $(\lambda, u_\lambda) $ such that $u_\lambda>0$ and $\|u_\lambda\|_{L^2(B_R)}^2=b$;
\item[$(ii)$] for any $0<c<b$, there exists at least two normalized solutions $(\lambda,u_\lambda),(\tilde{\lambda},u_{\tilde{\lambda}}) $, 
such that $u_\lambda>0$, $u_{\tilde{\lambda}}>0$ and 
$\|u_\lambda\|_{L^2(B_R)}^2=\|u_{\tilde{\lambda}}\|_{L^2(B_R)}^2=c$;
\item[$(iii)$]for any $c>b$, there exists no normalized solution $(\lambda,u_{\lambda})$ such that  
     $u_{\lambda}>0$ and $\|u_\lambda\|_{L^2(B_R)}^2=c$.
  \end{itemize}

  \item[$(2)$]
  If $p=2+\frac{4}{N}$, then there exists $D \geq d=\|Q\|_{L^2(\mathbb R^N)}>0$ such that 
\begin{itemize}
  \item[$(i)$] for any $0<c < d$, there exists at least one normalized solution $(\lambda,u_\lambda)$, 
  such that $u_\lambda>0$ and $\|u_\lambda\|_{L^2(B_R)}^2=c$;
  \item[$(ii)$] for any $c>D$, there exists no normalized solution $(\lambda,u_{\lambda})$, 
  such that $u_\lambda>0$ and $\|u_\lambda\|_{L^2(B_R)}^2=c$;
\end{itemize} 
where $Q$ is a ground state of \eqref{global}.

  \item[$(3)$]
  If $2<p<2+\frac{4}{N}$, then for any $0<c<+\infty$, there exists at least one normalized solution $(\lambda,u_\lambda) $, 
  such that $u_\lambda>0$ and $\|u_\lambda\|_{L^2(B_R)}^2=c$.

\end{itemize} 

\end{theorem}

Consider the equation 
\begin{align}\label{1plusx2}
  \begin{cases}-\Delta u +\lambda u= \frac{1}{(1+|x|^k)^s} u^p \quad \text{ in }B,
    \\ u>0, \quad u_{|\partial B}=0, \end{cases}
\end{align}
where $B$ is the unit ball, $2<p<2^*$, $s>1$ and $k\geq 0$.
\cite{GidasB} yields the radial symmetry and monotonicity of positive solutions of \eqref{1plusx2}, and it is not difficult to 
check that $(\Omega)$ and $(f_1)$-$(f_3)$ hold.

\begin{theorem}

\begin{itemize}

\item[$(1)$]Assume $N=3,-\lambda_1<\lambda<0, 2<p\leq4$ $and$ $ks\leq4-p$, then \eqref{1plusx2} admits a unique positive solution, 
which is radially symmetric and decreasing with respect to the radial variable $r=|x|.$

\item[$(2)$]  Assume $N=2,-\lambda_1<\lambda<0, 2<p\leq6$ $and$ $2ks\leq6-p$, then \eqref{1plusx2} admits a unique positive solution, 
which is radially symmetric and decreasing with respect to $the$ radial $variable\:r=|x|.$

\item[$(3)$]  Assume $N\geq2$ and $\lambda\geq0$, then \eqref{1plusx2} admits a unique positive solution, which is radially symmetric 
and decreasing with respect to the radial variable $r=|x|$.

\end{itemize} 
\end{theorem}

\begin{theorem}

\begin{itemize}
\item[$(1)$] Assume that
$$
N=2,4<p\leq6~~\text{and}~~2ks\leq6-p
$$
or
$$
N=3,\frac{10}3<p\leq4~~\text{and}~~ks\leq4-p,
$$
then
\begin{itemize}
  \item[$(i)$] \eqref{1plusx2} admits at least one positive normalized solution $(\lambda, u_\lambda)$ such that $\|u_\lambda\|_{L^2(B)}^2=b$;
  \item[$(ii)$] for any $0<c<b$, \eqref{1plusx2} admits at least two positive normalized solutions 
  $(\lambda,u_\lambda),(\tilde{\lambda},u_{\tilde{\lambda}})$, such that $\|u_\lambda\|_{L^2(B)}^2=\|u_{\tilde{\lambda}}|\|_{L^2(B)}^2=c$;
  \item[$(iii)$]for any $c>b$, \eqref{1plusx2} admits no positive normalized solution $(\lambda,u_{\lambda})$ 
  such that $\|u_\lambda\|_{L^2(B)}^2=c$, 
      where $u_{\lambda}>0$ and $\|u_\lambda\|_{L^2(B)}^2=c$.
    \end{itemize}
    \item[$(2)$]Assume that
    $$
    N=2,p=4~~~\text{and}~~~2ks\leq2, 
    $$
    or
    $$
    N=3,p=\frac{10}{3}~~~\text{and}~~~ks\leq\frac{2}{3},
    $$
    
    then there exists $D \geq d=\|Q\|_{L^2(\mathbb R^N)}>0$ such that 
 \begin{itemize}
    \item[$(i)$] for any $0<c < d$, \eqref{1plusx2} admits at least one normalized solution $(\lambda,u_\lambda) $, 
    such that $\|u_\lambda\|_{L^2(B)}^2=c$;
    \item[$(ii)$] for any $c>D$, there exists no normalized solution $(\lambda,u_{\lambda})$, such that $\|u_\lambda\|_{L^2(B)}^2=c$;
  \end{itemize} 
  where $Q$ is the unique positive solution of \eqref{global}.
  
    \item[$(3)$]Assume that
    $$
    N=2,2<p<4~~~\text{and}~~~2ks\leq6-p
    $$
    or
    $$
    N=3,2<p<\frac{10}{3}~~~\text{and}~~~ks\leq4-p,
    $$
     then for any $0<c<+\infty$, there exists at least one normalized solution $(\lambda,u_\lambda) $, 
    such that $\|u_\lambda\|_{L^2(B)}^2=c$.

  \end{itemize}

\end{theorem}

\medskip

\section{Orbital stability and instability}
In this section, we focus on the orbital stability and instability of the standing wave solutions 
$e^{i\lambda t}u_\lambda(x)$ for \eqref{c:intro}.

Recall that the solutions are called orbitally stable if for each $\epsilon>0$ there exists $\delta>0$ such that, 
whenever $\phi_0\in H_0^1(\Omega,\mathbb{C})$ is such that 
$\|\phi_0- u_\lambda\|_{H_0^1(\Omega,\mathbb{C}) }< \delta$ and $\phi(t,x)$ is the solution of ~\eqref{orbit}
with $\phi(0,\cdot)=\phi_0$ in some interval $[0,t_0)$, then $\phi(t,\cdot)$ can be continued to a solution in $0\leq t<\infty$ and

$$
\sup_{0<t<\infty}\inf_{w\in\mathbb{R}}\|\phi(t,x)-e^{i\lambda w}u_\lambda(x)\|_{H_0^1(\Omega,\mathbb{C})}<\epsilon;
$$
otherwise, they are called unstable.

\begin{lemma}\label{curve}
Assume that $(\Omega)$, $(f_1)-(f_2)$ and (N) hold.
  Set $K=\{(\lambda,u):\lambda> -\lambda_1, u \in K_\lambda, u>0\}$.
  If $K$ has only one element for any $\lambda>-\lambda_1$,
then $K$ is a $C^1$ curve in $\mathbb R \times H_0^1(\Omega)$.

\end{lemma}

\begin{proof}
  Similar to the proof of~\cite[Theorem 18]{Shatah}, by implicit function theorem, 
  we can prove that $K$ is a $C^1$ curve in $\mathbb R \times H_0^1(\Omega)$. 

\end{proof}

To study the orbital stability, we introduce the following result, in the form of the abstract theory developed in\cite{Grillakis}:

\begin{proposition}\label{abstracttheory}

Under the hypotheses of Lemma ~\ref{curve} and assume that (LWP) holds. Set $m(\lambda)=\|u_\lambda\|_{L^2(\Omega)}$, 
Then if  $m'(\lambda) > 0$ (respectively $< 0$), the standing wave
$e^{i\lambda t} u_\lambda(x)$ is orbitally stable (respectively unstable).

\end{proposition}

\begin{proof} [\textbf{Proof of Theorem~\ref{orbitalstability}}]

  For $2+\frac{4}{N}<p<2^*$, take $\lambda^*\in(-\lambda_1, +\infty)$ such that 
  $$
  m( \lambda^* ) = \max _{\lambda\in (-\lambda_1, +\infty) }m(\lambda),
  $$
   and $c\in \lambda\in(-\infty,\lambda_1)$
   $(0,m(\lambda^*))$ is a regular value of $m$, i.e. $m^{-1}(c)=\{\lambda_c^{(1)},\cdots,\lambda_c^{(k)}\}$ 
   and $m^{\prime}(\lambda_c^{(i)})\neq0,i=$ $1,2,\cdots,k.$ By Sard's theorem, regular values of $m$ are almost everywhere 
   in $(0,m(\lambda^*))$. 
   Since
   $$
   \lim_{\lambda\to-\lambda_1}m(\lambda)=\lim_{\lambda\to+\infty}m(\lambda)=0,
   $$
   $k$ must be even and we set $k=2s$. 
   Assume $\lambda_c^{(1)}<\lambda_c^{(2)}<\cdots<\lambda_c^{(2s)}$ .
   Then $m(\lambda)<m(\lambda_c^{(1)})$ if $-\lambda_1<\lambda< \lambda_c^{( 1) }$ and thus $m^{\prime}(\lambda_c^{(1)})>0.$ 
   Similarly, $m( \lambda) < m( \lambda_c^{( 2s) }) $ if $ \lambda >\lambda_c^{( 2s) }$ and 
   thus $m^{\prime}(\lambda_c^{(2s)})<0.$ By proposition ~\ref{abstracttheory}, the standing wave 
   $e^{i\lambda_c^{(1)}t}u_{\lambda_c^{(1)}(x)}$ is orbitally stable while  $e^{i\lambda_c^{(2s)}t}u_{\lambda_c^{(2s)}}(x)$
   is orbitally unstable.
   
   For $2<p \leq 2+\frac{4}{N}$, similar to the proof above, we can obtain the orbital stability results as well.

 \end{proof} 

 \medskip

\textbf{Acknowledgements}
I would like to thank C. Li and S. J. Li for fruitful disscussions and constant support during the development of this work.

\bigskip


\end{document}